\documentclass[10pt]{amsart}
\usepackage{latexsym,amsmath,amsfonts,amscd,amssymb}
\usepackage{graphics}
\newtheorem{theorem}{Theorem}[section]

\theoremstyle{remark}

\newcommand{\cC}{{\mathcal C}}
\newcommand{\cD}{{\mathcal D}}

\newcommand{\CC}{{\mathbb C}}

\title{A simple proof of the Fundamental Theorem of Algebra}


\author[R. P\'{e}rez-Marco]{Ricardo P\'{e}rez-Marco}
\address{CNRS, IMJ-PRG, Unversit\'e de Paris, Bo\^\i te courrier 7012, 75005 Paris Cedex 13, France}
\email{ricardo.perez.marco@gmail.com}

\begin{document}

\begin{abstract}
We present a simple short proof of the Fundamental Theorem of Algebra, without complex analysis and with a minimal use of topology. It can be taught in a first year calculus class.
\end{abstract}

\maketitle

\section{Statement.}

\begin{theorem}
 A non constant polynomial $P(z)\in \CC[z]$ with complex coefficients has a root.
\end{theorem}

The proof is based only on the following elementary facts:

$\bullet$ A polynomial has at most a finite number of roots.

$\bullet$ The Implicit Function Theorem.

$\bullet$ Removing from $\CC$ a finite number of points leaves an open connected space.

\section{The proof.}
 
It is enough to consider a monic polynomial $P$. We denote by $\cC=(P')^{-1} (0)$ the finite set of critical points of $P$, and by 
$\cD=P(\cC)$ the finite set of critical values of $P$.

$\bullet$ Let  $R=\{ c\in \CC; \text{ the polynomial }P(z)-c \text{ has at least a simple root and no double roots} \}$.

$\bullet$ $R\subset \CC-\cD$. This is because if $c\in \cD$, then 
$c=P(z_0)$ for some critical point $z_0 \in \cC$, hence $P'(z_0)=0$ and $P(z)-c=0$ has a double root at $z_0$. Note that $\CC-\cD$ is open and connected ($\cD$ being finite).

$\bullet$ $R$ is open. This is an application of the Implicit Function Theorem. 
Let $c_0\in R \subset \CC-\cD$, and  $z_0\in \CC$ be a root of $P(z)-c_0$. We apply the Implicit Function Theorem to the 
equation $F(z,c)=P(z)-c=0$. Since $\frac{\partial F}{\partial z}(z_0,c_0)=P'(z_0)\not=0$, there is a neighborhood $U$ of $c_0$ such that for $c\in U$ we have a root $z(c)$ of 
$P(z)-c$. Taking $U$ small enough, by continuity of $P'$ and $c\mapsto z(c)$, we have $P'(z(c))\not=0$ and the root $z(c)$ is simple. 
Since $\CC-\cD$ is open we can take $U\subset \CC-\cD$ and $P(z)-c$ does not have 
any double root, thus $U\subset R$.

$\bullet$ $R$ is closed in $\CC-\cD$. Because $P$ is monic, if $c$ is uniformly bounded then any root of $P(z)-c$ is uniformly bounded (since $P(z)/z^n\to 1$ uniformly when $z\to \infty$,
if $n$ is the degree). 
We can take a subsequence of $c_n\to c_\infty\ \in \CC-\cD$ and a converging subsequence of  roots of $P(z)-c_n$. By continuity, the limit is a root of $P(z)-c_\infty$,  
so this polynomial has roots. 
Moreover, all roots of $P(z)-c_\infty$ are simple since  $c_\infty \in \CC-\cD$.

$\bullet$ $R$ is non-empty. For any $a\in \CC$ we have that for $c=P(a)$, $P(z)-c$ has at least $z=a$ as root. If we choose $a\in \CC-P^{-1}(\cD)$, then for 
any root $z_0$ of $P(z)-c$ with $c=P(a)$, we have $P(z_0)=P(a) \notin \cD$, so $z_0\notin P^{-1}(\cD)$, but $\cC\subset P^{-1}(\cD)$, and $z_0\notin \cC$, and the root $z_0$ 
is simple.

The above proves that $R=\CC-\cD$. Now, if $0\in \cD$, then $0=P(z_0)$ for a critical point $z_0$ of $P$ that is also a root of $P$. If $0\notin \cD$, 
then $0\in R=\CC-\cD$ and the equation $P(z)-0=0$ has a simple root. In all cases $P$ has a root. $\diamond$

\section{Comment.}

The above proof is inspired from a beautiful proof by Daniel Litt \cite{Li}. He works in the global space of monic polynomials of degree $n\geq 1$ 
(biholomorphic to $\CC^n$), and removes the algebraic locus $\cD_n$, defined by the discriminant, of 
polynomials with a double root. He uses that the complement of an algebraic variety in $\CC^n$ is connected. Essentially the proof above 
achieves the same goal in a more elementary way working with $n=1$. In particular, we only need the simpler fact that the complement 
of a finite set in the plane is connected (which for $n=1$ is the 
same as the connectedness of the complement of an algebraic variety in $\CC^n$). We also avoid the use of discriminants.

\bigskip

\textbf{Acknowledgment.} I am grateful to my friends Marie-Claude Arnaud, Kingshook Biswas, Alain Chenciner and Yann Levagnini for their comments and suggestions to improve the presentation.  
In particular, to Kingshook that proposed a simplification of a first draft.


\begin{thebibliography}{9}



\bibitem{Li} LITT, D.; {\it Yet another proof of the Fundamental Theorem of Algebra}, Manuscript, 2011.

(www.daniellitt.com/blog/2016/10/6/a-minimal-proof-of-the-fundamental-theorem-of-algebra)




\end{thebibliography}
\end{document}